\documentclass[12pt,reqno]{article}
\usepackage{amssymb,amscd,amsmath,amsthm,color}
\newtheorem{theorem}{Theorem}[section]
\newtheorem{lemma}[theorem]{Lemma}
\newtheorem{cor}[theorem]{Corollary}
\newtheorem{prop}[theorem]{Proposition}
\newtheorem{conjecture}[theorem]{Conjecture}
\newtheorem{question}[theorem]{Question}
\usepackage{enumerate,amssymb}
\theoremstyle{definition}
\newtheorem{definition}[theorem]{Definition}
\newtheorem{remark}[theorem]{Remark}

\theoremstyle{remark}

\numberwithin{equation}{section}

\def\bC{\mathbb{C}}
\def\bE{\mathbb{E}}
\def\bM{\mathbb{M}}

\begin{document}
\baselineskip=15pt

\title{ Pinchings and  positive linear maps  }

\author{Jean-Christophe Bourin{\footnote{Supported by ANR 2011-BS01-008-01.}}\, and  Eun-Young Lee{\footnote{Research  supported by  Basic Science Research Program
through the National Research Foundation of Korea (NRF) funded by the Ministry of Education,
Science and Technology (2010-0003520)}}}

\date{ }

\maketitle

\vskip 10pt\noindent
{\bf Abstract.}  We employ  the pinching theorem, ensuring that some   operators $A$ admit {\it any} sequence of contractions as an operator diagonal of $A$,  to deduce/improve  two recent theorems of Kennedy-Skoufranis and Loreaux-Weiss for conditional expectations onto a masa  in the algebra of   operators on a Hilbert space. We also get a few results for sums in a unitary orbit.


\section{The pinching theorem}

We recall two theorems which are fundamental in the next sections to obtain several results  about positive linear maps, in particular conditional expectations, and unitary orbits.
These theorems were established in \cite{B}, we also
  refer to this article for various definitions and properties of the essential numerical range $W_e(A)$ of an operator $A$ in the algebra ${\mathrm{L}}({\mathcal{H}})$ of all
 (bounded linear) operators on an infinite dimensional, separable (real or complex) Hilbert space ${\mathcal{H}}$. 

We denote by ${\mathcal{D}}$ the unit disc of $\bC$. We write $A\simeq B$ to mean that the operators $A$ and $B$ are unitarily equivalent. This relation is extended to operators possibly acting on different Hilbert spaces, typically, $A$ acts on ${\mathcal{H}}$ and $B$ acts on an infinite dimensional subspace ${\mathcal{S}}$ of ${\mathcal{H}}$, or on the spaces ${\mathcal{H}}\oplus {\mathcal{H}}$ or $\oplus^{\infty}{\mathcal{H}}$.

\vskip 5pt
\begin{theorem}\label{pinching} Let $A\in{\mathrm{L}}({\mathcal{H}}) $ with $W_e(A)\supset{\mathcal{D}}$ and  $\{X_i\}_{i=1}^{\infty}$  a sequence in ${\mathrm{L}}({\mathcal{H}})$ such that
$\sup_i\|X_i\|<1$. Then,  a  decomposition  ${\mathcal{H}}=\oplus_{i=1}^{\infty}{\mathcal{H}}_i$ holds with 
 $A_{{\mathcal{H}}_i}\simeq X_i$ for all $i$.
\end{theorem}

\vskip 5pt
Of course, the direct sum refers to an orthogonal decomposition, and $A_{{\mathcal{H}}_i}$ stands for the compression of $A$ onto the subspace ${\mathcal{H}}_i$.

Theorem \ref{pinching}  tells us that we have a unitary congruence between an operator in ${\mathrm{L}}(\oplus^{\infty}{\mathcal{H}})$ and a "pinching" of $A$,
$$
\bigoplus_{i=1}^{\infty} X_i \simeq \sum_{i=1}^{\infty} E_iAE_i
$$
for some sequence of mutually orthogonal infinite dimensional projections $\{E_i\}_{i=1}^{\infty}$ in ${\mathrm{L}}({\mathcal{H}})$ summing up to the identity $I$. Thus   $\{X_i\}_{i=1}^{\infty}$ can be regarded as an operator diagonal of $A$. In particular, if $X$ is an operator on ${\mathcal{H}}$ with $\| X\| <1$, then, $A$ is unitarily congruent to an operator on ${\mathcal{H}}\oplus {\mathcal{H}}$ of the form,
\begin{equation}\label{compression}
A\simeq \begin{pmatrix} X&\ast \\ \ast& \ast\end{pmatrix}.
\end{equation}

For a sequence of normal operators, Theorem \ref{pinching} admits a variation. Given   ${\mathcal{A}},{\mathcal{B}}\subset \bC$, the notation  ${\mathcal{A}}\subset_{st}{\mathcal{B}}$ means that ${\mathcal{A}}+r{\mathcal{D}}\subset {\mathcal{B}}$ for some $r>0$.

\vskip 5pt
\begin{theorem}\label{pinchingnormal} Let $A\in{\mathrm{L}}({\mathcal{H}}) $ with $W_e(A)\supset{\mathcal{D}}$ and  $\{X_i\}_{i=1}^{\infty}$  a sequence of normal operators in ${\mathrm{L}}({\mathcal{H}})$ such that
$\cup_{i=1}^{\infty}W(X_i)\subset_{st} W_e(A)$. Then,  a decomposition  ${\mathcal{H}}=\oplus_{i=1}^{\infty}{\mathcal{H}}_i$ holds with  
 $A_{{\mathcal{H}}_i}\simeq X_i$ for all $i$. 
\end{theorem}

\vskip 5pt 
In Section 3, our concern is the study of generalized diagonals, i.e., conditional expectations onto a masa in ${\mathrm{L}}({\mathcal{H}}) $, of the unitary orbit of an operator. The pinching theorems are the good tools for this study; we easily  obtain and considerably improve two recent theorems, of Kennedy and Skoufranis 
for normal operators, and Loreaux and Weiss for idempotent operators. 
Section 4 deals with an application to the class of unital, positive linear maps which are trace preserving. 
Section, 5 collects a few questions on possible extension of Theorems \ref{pinching} and
\ref{pinchingnormal} in the setting of von Neumann algebras.

The next section gives applications which only require  \eqref{compression}.  These results mainly focus on sums of two operators in a unitary orbit. 

\section{Sums in a unitary orbit}

\vskip 5pt\noindent
We recall  a straightforward consequence of \eqref{compression} for the weak convergence, \cite[Corollary 2.4]{B}.

\vskip 5pt
\begin{cor}\label{cor1seq} Let $A,X\in{\mathrm{L}}({\mathcal{H}})$ with  $W_e(A)\supset{\mathcal{D}}$ and $\| X\| \le 1$. Then there exists a sequence of unitaries $\{U_n\}_{n=1}^{\infty}$  in ${\mathrm{L}}({\mathcal{H}})$ such that 
$${\mathrm{wot}}\!\!\! \lim_{n\to+\infty}U_nAU_n^*= X. $$
\end{cor} 

\vskip 5pt
Of course, we cannot replace the weak convergence by the strong convergence; for instance if  $A$ is invertible and $\| Xh\| <\| A^{-1}\|^{-1}$ for some unit vector $h$, then $X$ cannot be a strong limit from the unitary orbit of $A$. However, the next best thing does happen. Moreover, this is even true for the $\ast$-strong operator topology.

\vskip 5pt
\begin{cor}\label{cor2seq} Let $A,X\in{\mathrm{L}}({\mathcal{H}})$ with  $W_e(A)\supset{\mathcal{D}}$ and $\| X\| \le 1$. Then there exist two sequences of unitaries $\{U_n\}_{n=1}^{\infty}$ and $\{V_n\}_{n=1}^{\infty}$ in ${\mathrm{L}}({\mathcal{H}})$ such that 
$$\ast\,{\mathrm{sot}}\!\!\! \lim_{n\to+\infty}\frac{U_nAU_n^*+V_nAV_n^* }{2}= X. $$
\end{cor} 

\vskip 5pt
\begin{proof} From \eqref{compression} we also have
\begin{equation*}
A\simeq \begin{pmatrix} X&-R \\ -S&T\end{pmatrix}.
\end{equation*}
Hence there exist two unitaries $U,V:{\mathcal{H}}\to{\mathcal{H}}\oplus{\mathcal{H}}$
such that 
\begin{equation}\label{eq1}
\frac{UAU^*+VAV^*}{2}= \begin{pmatrix} X&0 \\ 0&T\end{pmatrix}.
\end{equation}
Now let $\{e_n\}_{n=1}^{\infty}$ be a basis of  ${\mathcal{H}}$ and choose any unitary $W_n:{\mathcal{H}}\oplus{\mathcal{H}}\to {\mathcal{H}} $ such that $W_n(e_j\oplus 0)=e_j$ for all $ j\le n$. Then
$$
X_n:=W_n  \begin{pmatrix} X&0 \\ 0&T\end{pmatrix} W_n^*
$$
strongly converges to $X$. Indeed, $\{X_n\}$ is bounded in norm and, for all $j$, $X_ne_j\to Xe_j$. Taking adjoints,
$$
X_n^*=W_n  \begin{pmatrix} X^*&0 \\ 0&T^*\end{pmatrix} W_n^*,
$$
we also have $X_n^*\to X$ strongly. Setting $U_n=W_nU$ and $V_n=W_nV$ and using \eqref{eq1} completes the proof. \end{proof}

\vskip 5pt 
\begin{remark}\label{remmean} Corollary \ref{cor2seq} does not hold for the convergence in norm. We give an example. Consider the permutation matrix
$$
T=\begin{pmatrix} 0&0&1 \\ 1&0&0 \\ 0&1&0 \end{pmatrix}
$$
and set $A=2\oplus^{\infty} T$ regarded as an operator in ${\mathrm{L}}({\mathcal{H}})$. Then  $W_e(A)\supset{\mathcal{D}}$, however $X=(1/2)I$ is not a norm limit from the unitary orbit of $A$. Equivalently, $(1/2)I$ is not a norm limit from the unitary orbit of $(A+A^*)/2$. Indeed,  $(A+A^*)/2= I-(3/2)P$ for some projection $P$.
\end{remark}

\vskip 5pt 
\begin{remark} The converse of Corollary \ref{cor2seq} holds: if $A\in {\mathrm{L}}({\mathcal{H}})$ has the property  that any contraction is a strong limit of a mean of two operators in its unitary orbit, then necessarily $W_e(A)\supset {\mathcal{D}}$. This is checked by arguing as in the proof of Corollary \ref{cordoubly}.
\end{remark}

\vskip 5pt
We reserve the word "projection" for selfadjoint idempotent. A strong limit of idempotent operators is still idempotent; thus, the next corollary is rather surprising.

\vskip 5pt
\begin{cor}\label{idem2seq} Fix $\alpha>0$. There exists an idempotent $Q\in{\mathrm{L}}({\mathcal{H}})$ such that for every $X\in{\mathrm{L}}({\mathcal{H}})$ with $\| X\| \le \alpha$ we have
 two sequences of unitaries $\{U_n\}_{n=1}^{\infty}$ and $\{V_n\}_{n=1}^{\infty}$ in ${\mathrm{L}}({\mathcal{H}})$ for which
$$\ast\,{\mathrm{sot}}\!\!\! \lim_{n\to+\infty} U_nQU_n^*+V_nQV_n^* = X. $$
\end{cor}

\vskip 5pt
\begin{proof}
Let $a>0$, define a two-by-two idempotent matrix
\begin{equation}\label{eqidempotent}
M_a=\begin{pmatrix} 1&0 \\ a&0 \end{pmatrix}
\end{equation}
and set $Q=\oplus^{\infty}M_a$  regarded as an operator in ${\mathrm{L}}({\mathcal{H}})$. Since the numerical range $W(\cdot)$ of
$$
\begin{pmatrix} 0&0 \\ 2&0 \end{pmatrix}
$$
is ${\mathcal{D}}$, we infer that $W(2\alpha^{-1}M_a)=W_e((2\alpha^{-1}Q)\supset{\mathcal{D}}$ for a large enough $a$. The result then follows from Corollary \ref{cor2seq} with $A=2\alpha^{-1}Q$  and the contraction $\alpha^{-1}X$. 
 \end{proof}

\vskip 5pt 
 Corollary \ref{idem2seq} does not hold for the convergence in norm.

\vskip 5pt
\begin{prop}\label{prop1} Let $X\in{\mathrm{L}}({\mathcal{H}})$ be of the form $\lambda I +K$ for a compact operator $K$
  and a scalar $\lambda\notin \{0,1,2\}$. Then $X$ is not  norm limit of $U_nQU_n^*+V_nQV_n^*$ for any sequences of unitaries $\{U_n\}_{n=1}^{\infty}$ and $\{V_n\}_{n=1}^{\infty}$ and any idempotent $Q$ in ${\mathrm{L}}({\mathcal{H}})$. 
\end{prop}

\vskip 5pt
\begin{proof} First observe that if $\{A_n\}_{n=1}^{\infty}$ and $\{B_n\}_{n=1}^{\infty}$ are two bounded sequences in ${\mathrm{L}}({\mathcal{H}})$ such that $A_n-B_n\to 0$ in norm, then we also have $A_n^2-B_n^2\to 0$ in norm; indeed 
$$
A_n^2-B_n^2= A_n(A_n-B_n) + (A_n-B_n)B_n.
$$
Now, suppose that $\lambda\neq 1$ and that we have  the (norm) convergence,
$$
U_nQU_n^*+V_nQV_n^*\to \lambda I +K.
$$
Then we also have 
\begin{equation}\label{eqL1}
W_nQW_n^*-\left(-Q+\lambda I +U_n^*KU_n\right)\to 0
\end{equation}
where $W_n:=U_n^*V_n$. Hence,
 by the previous observation,
$$
(W_nQW_n^*)^2-\left(-Q+\lambda I +U_n^*KU_n\right)^2\to 0,
$$
that is
\begin{equation}\label{eqL2}
W_nQW_n^*-\left(-Q+\lambda I +U_n^*KU_n\right)^2\to 0.
\end{equation}
Combining \eqref{eqL1} and \eqref{eqL2} we get
$$
\left(-Q+\lambda I +U_n^*KU_n\right)-\left(-Q+\lambda I +U_n^*KU_n\right)^2\to 0
$$
hence
$$
(-2+2\lambda)Q +(\lambda -\lambda^2 )I + K_n\to 0
$$
for some bounded sequence of compact operators $K_n$. Since $\lambda\neq 1$, we have
$$
Q=\frac{\lambda}{2} I+L
$$
for some compact operator $L$. Since $Q$ is idempotent, either $\lambda=2$ or $\lambda=0$.
\end{proof}

\vskip 5pt The operator $X$ in Proposition \ref{prop1} has the special property that $W_e(X)$ is reduced to a single point. However Proposition \ref{prop1} may also hold when $W_e(X)$ has positive measure.

\vskip 5pt
\begin{cor} Let $Q$ be an  idempotent in ${\mathrm{L}}({\mathcal{H}})$ and $z\in\bC\setminus\{0,1,2\}$. Then, there exists $\alpha >0$ such that the following property holds:  
\begin{itemize}
\item[]
If $X\in{\mathrm{L}}({\mathcal{H}})$ satisfies $\| X-zI\| \le \alpha$, then $X$ is not norm limit of $U_nQU_n^*+V_nQV_n^*$ for any
  sequences of unitaries $\{U_n\}_{n=1}^{\infty}$ and $\{V_n\}_{n=1}^{\infty}$ in ${\mathrm{L}}({\mathcal{H}})$.
\end{itemize}
\end{cor}

\vskip 5pt
\begin{proof} By the contrary, $zI$ would be a norm limit of $U_nQU_n^*+V_nQV_n^*$ for some unitaries $U_n,V_n$, contradicting Proposition \ref{prop1}.
\end{proof}

\vskip 5pt More operators with large numerical and essential numerical ranges are given in the next proposition. An operator $X$ is stable when its real part $(X+X^*)/2$  is negative definite (invertible).

\vskip 5pt
\begin{prop}\label{prop2} If $X\in{\mathrm{L}}({\mathcal{H}})$ is stable, then $X$ is not  norm limit of $U_nQU_n^*+V_nQV_n^*$ for any sequences of unitaries $\{U_n\}_{n=1}^{\infty}$ and $\{V_n\}_{n=1}^{\infty}$ and any idempotent $Q$ in ${\mathrm{L}}({\mathcal{H}})$. 
\end{prop}

\vskip 5pt
\begin{proof} We have a decomposition ${\mathcal{H}}={\mathcal{H}}_s\oplus {\mathcal{H}}_{ns}$ in two invariant subspaces of $Q$ such that $Q$ acts on ${\mathcal{H}}_s$ as a selfadjoint projection $P$, and $Q$ acts on ${\mathcal{H}}_{ns}$
as a purely nonselfadjoint idempotent, that is $A_{{\mathcal{H}}_{ns}}$ is unitarily equivalent to an operator on ${\mathcal{F}}\oplus{\mathcal{F}}$ of the form
\begin{equation}\label{purelyns}
Q_{{\mathcal{H}}_{ns}}\simeq\begin{pmatrix} I& 0 \\ R&0\end{pmatrix}
\end{equation}
where $R$ is a nonsingular positive operator  on a Hilbert space ${\mathcal{F}}$, so 
\begin{equation}\label{purelyns2}
Q\simeq P\oplus \begin{pmatrix} I& 0 \\ R&0\end{pmatrix}.
\end{equation}

Let $Y$ be  a norm limit of the sum of two sequences in the unitary orbit of $Q$. 
If the purely non-selfadjoint part ${\mathcal{H}}_{ns}$ is vacuous, then $Y$ is positive, hence $Y\neq X$.  If ${\mathcal{H}}_{ns}$ is not vacuous, \eqref{purelyns2} shows that
\begin{align*}
Q+Q^*&\simeq  2P\oplus \begin{pmatrix} 2I& R \\ R&0\end{pmatrix} \\
&\simeq  2P\oplus\left\{ \begin{pmatrix} I& I \\ I&I\end{pmatrix}
+  \begin{pmatrix} R& 0 \\ 0&-R\end{pmatrix}\right\}.
\end{align*} 
This implies that $\| (Q+Q^*)_+\| \ge \| (Q+Q^*)_-\|$, therefore $Y+Y^*$ cannot be negative definite, hence $X\neq Y$.
\end{proof} 
 
\vskip 5pt 
 It is known \cite{PT} that any operator is the sum of five idempotents.
We close this section by asking whether Corollorary \ref{idem2seq} admits a substitute for Banach space operators.

\vskip 5pt
\begin{question} Let ${\mathcal{ X}}$ be a separable Banach space and $T\in{\mathrm{L}}({\mathcal{ X}})$, the linear operators on ${\mathcal{ X}}$. Do there exist two sequences $\{P_n\}_{n=1}^{\infty}$ and $\{Q_n\}_{n=1}^{\infty}$ of idempotents in ${\mathrm{L}}({\mathcal{ X}})$ such that $T={\mathrm{sot}} \lim_{n\to+\infty} (P_n +Q_n)$ ?
\end{question}

\section{Conditional expectation onto a masa}

\subsection{Conditional expectation of general operators}

Kennedy and Skoufranis have studied the following problem: Let ${\mathfrak{X}}$ be a  maximal abelian $\ast$-subalgebra (masa) of a von Neumann algebra ${\mathfrak{M}}$, with corresponding expectation $\bE_{\mathfrak{X}}: {\mathfrak{M}}\to  {\mathfrak{X}}$ (i.e., a unital positive linear map such that $\bE_{\mathfrak{X}}(XM)= X\bE_{\mathfrak{X}}(M)$ for all $X\in{\mathfrak{X}}$ and $M\in{\mathfrak{M}}$). Given a normal operator $A\in{\frak{M}}$, determine the image by $\bE_{\frak{X}}$ of the unitary orbit of $A$, 
$$
\Delta_{\mathfrak{X}}(A)=\{\, \bE_{\mathfrak{X}}(UAU^*)\  : \ U \, {\mathrm{a\ unitary\ in}}\ {\mathfrak{M}}\, \}.
$$
In several cases, they  determined the norm closure of $\Delta_{\frak{X}}(A)$. In particular,  \cite[Theorem 1.2]{KS} can be stated in the following two propositions.

\vskip 5pt
\begin{prop}\label{propKS} Let ${\mathfrak{X}}$ be a  masa in ${\mathrm{L}}({\mathcal{H}})$,  $X\in {\mathfrak{X}}$, and  $A$  a normal operator in $ {\mathrm{L}}({\mathcal{H}})$.  If $\sigma(X)\subset {\mathrm{conv}}\sigma_e(A)$,  then $X$ lies in the norm closure of $\Delta_{\frak{X}}(A)$. 
\end{prop}

\vskip 5pt
\begin{prop}\label{propKS2} Let ${\mathfrak{X}}$ be a continuous masa in ${\mathrm{L}}({\mathcal{H}})$,  $X\in {\mathfrak{X}}$, and $A$  a normal operator in $ {\mathrm{L}}({\mathcal{H}})$. If $X$ lies in the norm closure of $\Delta_{\frak{X}}(A)$, then $\sigma(X)\subset {\mathrm{conv}}\sigma_e(A)$.
\end{prop}

\vskip 5pt
Since we deal with normal operators, $\sigma(X)\subset {\mathrm{conv}}\sigma_e(A)$ means $W(X)\subset W_e(A)$. 
 Proposition \ref{propKS2} needs the continuous assumption. It is a rather simple fact;   we generalize it in Lemma \ref{lemred}: {\it  Conditional expectations reduce  essential numerical ranges,   $W(\bE_{\frak{X}}(T))\subset W_e (T)$ for all $T\in {\mathrm{L}}({\mathcal{H}})$}. Thus, the main point of \cite[Theorem 1.2]{KS} is  Proposition \ref{propKS} which says that if $W(X)\subset W_e(A)$ then $X$ can be approximated by operators   of the form $\bE_{\frak{X}}(UAU^*)$ with unitaries $U$. With the slightly stronger assumption $W(X)\subset_{st} W_e(A)$, Theorem \ref{pinchingnormal} guarantees, via the following corollary, that $X$ is exactly  of this form, furthermore the normality assumption on $A$ is not necessary.

\vskip 10pt
\begin{cor}\label{corKS} Let ${\mathfrak{X}}$ be a  masa in ${\mathrm{L}}({\mathcal{H}})$,  $X\in {\mathfrak{X}}$ and $A\in {\mathrm{L}}({\mathcal{H}})$. 
If $W(X)\subset_{st} W_e(A)$, then $X=\bE_{\mathfrak{X}}(UAU^*)$ for some unitary operator $U\in{\mathrm{L}}({\mathcal{H}})$.
\end{cor} 

\vskip 5pt
\begin{proof} First, we note a simple fact: Let $\{P_i\}_{i=1}^{\infty}$ be a sequence of orthogonal projections in ${\mathfrak{X}}$ such that $\sum_{i=1}^{\infty}P_i=I$, and let $Z\in {\mathrm{L}}({\mathcal{H}})$ such that $P_iZP_i\in {\mathfrak{X}}$ for all $i$. Then, we have a strong sum
$$\bE_{\mathfrak{X}}(Z)=\sum_{i=1}^{\infty} P_iZP_i.$$

Now, denote by ${\mathcal{H}}_i$ the range of $P_i$ and assume $\dim{\mathcal{H}}_i=\infty$ for all $i$. We have  $W_e(X_{{\mathcal{H}}_i})\subset W_e(X)$, hence $$\cup_{i=1}^{\infty}W(X_{{\mathcal{H}}_i})\subset_{st} W_e(A).$$
We may then apply Theorem \ref{pinchingnormal} and get a unitary $U$ on ${\mathcal{H}}=\oplus_{i=1}^{\infty}{\mathcal{H}}_i$ such that
$$
A\simeq UAU^* =
\begin{pmatrix} 
X_{{\mathcal{H}}_1} &\ast & \cdots &\cdots  \\ 

\ast &X_{{\mathcal{H}}_2} &\ast &\cdots  \\ 

 \vdots &\ast  & \ddots &\ddots \\ 

\vdots  &\vdots  &\ddots  &\ddots \\ 
\end{pmatrix}.
$$
Since $0\oplus\cdots\oplus X_{{\mathcal{H}}_i}\oplus0\cdots \in {\mathfrak{X}}$ for all $i$, the previous simple fact shows that
$$
\bE_{\mathfrak{X}}(UAU^*)=\bigoplus_{i=1}^{\infty}X_{{\mathcal{H}}_i}= X.
$$
\end{proof}

\vskip 5pt
\begin{remark}
 Corollary \ref{corKS} also covers the assumption $W(X)\subset W_e(A)$ of Proposition \ref{propKS}. Indeed, $W_e(A)\subset_{st} W_e(A+D)$ for some  normal  operator $D$ with arbitrarily small norm, and we may apply Corollary \ref{corKS} to $X$ and $A+D$.
\end{remark}

\subsection{A reduction lemma}

The following result extends Proposition \ref{propKS2}, the "easy" part of Kennedy-Skoufranis' theorem \cite[Theorem 1.2]{KS}.

\vskip 5pt
\begin{lemma}\label{lemred} If ${\mathfrak{X}}$ is a  masa in ${\mathrm{L}}({\mathcal{H}})$ and $Z\in{\mathrm{L}}({\mathcal{H}})$, then $W(\bE_{\mathfrak{X}}(Z))\subset
\overline{W}(Z)$ and $W_e(\bE_{\mathfrak{X}}(Z))\subset
W_e(Z)$.
\end{lemma}

\vskip 5pt
\begin{proof} (1) Assume $Z$ is normal. We may identify the unital $C^*$-algebra ${\mathfrak{A}}$ spanned by $Z$ with $C^0(\sigma(Z))$ via a $\ast$-isomorphism $\varphi: C^0(\sigma(Z))\to {\mathfrak{A}}$ with $\varphi(z\mapsto z)=Z$. Let $h\in{\mathcal{H}}$ be a unit vector.
For $f\in C^0(\sigma(Z))$, set
$$
\psi(f)=\langle h, \bE_{\mathfrak{X}}(\varphi(f))h\rangle.
$$
Then $\psi$ is a positive linear functional on $C^0(\sigma(Z))$ and $\psi(1)=1$. Thus $\psi$ is a Radon measure induced by a  probabilty measure $\mu$,
$$
\psi(f)=\int_{\sigma(Z)} f(z)\, {\mathrm{d}}\mu(z)
$$
We then have $
\langle h, \bE_{\mathfrak{X}}(Z)h\rangle =\psi(z) \in {\mathrm{conv}}(\sigma(Z)).
$
Since ${\mathrm{conv}}(\sigma(Z))=\overline{W}(Z)$, we obtain $W(\bE_{\mathfrak{X}}(Z))\subset
\overline{W}(Z)$.

(2) Let $Z$ be a general operator in ${\mathrm{L}}({\mathcal{H}})$ and define a conditional expectation $$\bE_2:{\mathrm{L}}({\mathcal{H}}\oplus{\mathcal{H}})\to {\mathfrak{X}}\oplus{\mathfrak{X}}$$ by
$$
\bE_2\left( 
\begin{pmatrix} A&C \\ D& B\end{pmatrix}
\right)=
\begin{pmatrix} \bE_{\mathfrak{X}}(A)&0 \\ 0& \bE_{\mathfrak{X}}(B)\end{pmatrix}.
$$
From the first part of the proof, we infer
$$
W(\bE_{\mathfrak{X}}(Z))\subset W\left( 
\begin{pmatrix} (\bE_{\mathfrak{X}}( Z)&0 \\ 0& (\bE_{\mathfrak{X}}(B)\end{pmatrix}
\right)
\subset\overline{W}\left( 
\begin{pmatrix} Z&C \\ D& B\end{pmatrix}
\right)
$$
whener 
$
\begin{pmatrix} Z&C \\ D& B\end{pmatrix}
$
is normal. Since  we have, by a simple classical fact \cite{Hal},
$$
\overline{W}(Z)=\bigcap \overline{W}\left(\begin{pmatrix} Z&C \\ D& B\end{pmatrix}\right)
$$
where the intersection runs over all $B,C,D$ such that
$
\begin{pmatrix} Z&C \\ D& B\end{pmatrix}
$
is normal, we obtain $W(\bE_{\mathfrak{X}}(Z))\subset
\overline{W}(Z)$. 

(3) We deal with the essential numerical range inclusion. We can split ${\mathfrak{X}}$ into its discrete part ${\mathfrak{D}}$ and continuous part ${\mathfrak{C}}$ with the  corresponding decomposition of the Hilbert space,
$$
{\mathfrak{X}}={\mathfrak{D}}\oplus{\mathfrak{C}}, \qquad {\mathcal{H}}={\mathcal{H}}_d\oplus{\mathcal{H}}_c.
$$
We then have
\begin{equation}\label{red1}
W_e(\bE_{\mathfrak{X}}(Z))={\mathrm{conv}}\left\{
W_e(\bE_{\mathfrak{D}}(Z_{{\mathcal{H}}_d}))
;W_e(\bE_{\mathfrak{C}}(Z_{{\mathcal{H}}_c}))
 \right\}.
\end{equation}
We have an obvious inclusion
\begin{equation}\label{red2}
W_e(\bE_{\mathfrak{D}}(Z_{{\mathcal{H}}_d}))\subset W_e(Z_{{\mathcal{H}}_d}).
\end{equation}
On the other hand, for all compact operators $K\in{\mathrm{L}}({\mathcal{H}})$,
$$
W_e(\bE_{\mathfrak{C}}(Z_{{\mathcal{H}}_c}))=W_e(\bE_{\mathfrak{C}}(Z_{{\mathcal{H}}_c})+K_{{\mathcal{H}}_c})=W_e(\bE_{\mathfrak{C}}(Z_{{\mathcal{H}}_c}+K_{{\mathcal{H}}_c}))\subset\overline{W}(Z_{{\mathcal{H}}_c}+K_{{\mathcal{H}}_c})
$$
by the simple folklore fact that a conditional expectation onto a continous masa vanishes on compact operators  and part (2) of the proof. Thus, when $K$ runs over all compact operators, we obtain
\begin{equation}\label{red3}
W_e(\bE_{\mathfrak{D}}(Z_{{\mathcal{H}}_c}))\subset W_e(Z_{{\mathcal{H}}_c}).
\end{equation}
Combining \eqref{red1}, \eqref{red2} and  \eqref{red3} completes the proof.
\end{proof}

\subsection{Conditional expectation of idempotent operators}

\vskip 5pt
For discrete masas, unlike continuous masas \cite{KadSing}, there is a unique conditional expectation, which   merely consists in extracting the diagonal  with respect to an orthonormal basis. In a recent article,  Loreaux and Weiss give a detailed study of diagonals of idempotents in  ${\mathrm{L}}({\mathcal{H}})$. They established that a nonzero idempotent $Q$ has a zero diagonal with respect to some orthonormal basis if and only if $Q$ is not a Hilbert-Schmidt perturbation of a projection (i.e., a self-adjoint idempotent). They also showed that any sequence $\{a_n\}\in l^{\infty}$ such that $|a_n|\le \alpha$ for all $n$ and, for some $a_{n_0}$, $a_k=a_{n_0}$ for infinitely many $k$, one has a idempotent $Q$ such that $\| Q\|\le 18\alpha +4$ and $Q$ admits $\{a_n\}$ as a diagonal with respect to some orthonormal basis \cite[Proposition 3.4]{LW}. Using this, they proved that any sequence in $l^{\infty}$ is the diagonal of some idempotent operator \cite[Theorem 3.6]{LW}, answering a question of Jasper.  This statement is in the range of Theorem \ref{pinching}. Further, it is not necessary to confine to diagonals, i.e., discrete masas, and the constant $18\alpha +4$  can be improved; in the next corollary we explicit the best constant when $\alpha=1$.

\vskip 5pt
\begin{cor}\label{corLW} Let ${\mathfrak{X}}$   be a masa in ${\mathrm{L}}({\mathcal{H}})$ and $\alpha>0$. There exists an idempotent  $Q\in {\mathrm{L}}({\mathcal{H}})$,  such that for all  $X\in  {\mathfrak{X}}$ with $\| X\|<\alpha$,  we have  $X=\bE_{\frak{X}}(UQU^*)$ for some unitary operator $U\in{\mathrm{L}}({\mathcal{H}})$. If $\alpha=1$, $\| Q\|=\sqrt{5+2\sqrt{5}}$ is the smallest possible norm. 
\end{cor} 

\vskip 5pt
\begin{proof} As in the proof of Corollary \ref{idem2seq} we have an idempotent $Q$ such that $W_e(Q)\supset \alpha{\mathcal{D}}$,  hence the first and main part of  Corollary \ref{corLW} follows from
Corollary \ref{corKS}. The remaining parts require a few computations.

To obtain the bound $\sqrt{5+2\sqrt{5}}$ when $\alpha=1$ we get a closer look at
 $\oplus^{\infty} M_a$ with $M_a$ given by \eqref{eqidempotent} where $a$ is a positive scalar. We have
\begin{align*}
W(M_a)&=\left\{  \langle h, M_ah\rangle \ : \ h\in\bC^2, \|h\|=1\right\} \\
&=\left\{  |h_1|^2 + a\overline{h_2}h_1 \ : |h_1|^2 +|h_2|^2=1\right\}, 
\end{align*}
hence, with $h_1=re^{i\theta}$, $h_2=\sqrt{1-r^2}e^{i\alpha}$,
$$
W(M_a)=\bigcup_{0\le r\le 1}\left\{  r^2 + ar\sqrt{1-r^2}e^{i(\theta-\alpha)} : \ \theta, \alpha \in [0,2\pi] \right\}. 
$$
Therefore $W(M_a)$ is a union of circles $\Gamma_r$ with centers $r^2$ and radii $ar\sqrt{1-r^2}$. To have  ${\mathcal{D}}\subset W(M_a)$ it is necessary and sufficient  that $-1\in \Gamma_r$ for some $r\in [0,1]$, hence
\begin{equation}\label{min}
a=\frac{1+r^2}{r\sqrt{1-r^2}}.
\end{equation}
Now we minimize $a=a(r)$ given by \eqref{min} when $r\in(0,1)$ and thus obtain the matrix $M_{a_{\ast}}$ with smallest norm such that $W(M_{a_{\ast}})\supset{\mathcal{D}}$. Observe that $a(r)\to +\infty$ as $r\to 0$ and as $r\to 1$, and
$$
r^2(1-r^2)^{3/2}a'(r)= r^4 +4r^2-1.
$$
Thus $a(r)$  takes its minimal value $a_{\ast}$ when $r^2=\sqrt{5}-2$.
We have
$a_{\ast}^2=4+2\sqrt{5}$, hence 
$$\| M_{a_{\ast}}\|=\sqrt{5+2\sqrt{5}}.$$ 
Now, letting $Q=\oplus^{\infty} M_{a_{\ast}}$, we have  $W_e(Q)=W(M_{a_{\ast}})$, so that $Q$ is an idempotent in ${\mathrm{L}}({\mathcal{H}})$
such that  $W_e(Q)\supset {\mathcal{D}}$, and thus by Corollary \ref{corKS} any operator $X$  such that $\| X\|<1$ satifies $\bE_{\mathfrak{X}}(UQU^*)=X$ for some  unitary  $U$.

It remains to check that if $Q$ is an idempotent such that Corollary \ref{corLW} holds for any operator $X$ such that $\| X\|<1$, then $\| Q\| \ge \sqrt{5+2\sqrt{5}}$. To this end, we consider the purely nonselfadjoint part $Q_{{\mathcal{H}}_{ns}}$ of $Q$  in \eqref{purelyns},
\begin{equation*}
Q_{{\mathcal{H}}_{ns}}\simeq\begin{pmatrix} I& 0 \\ R&0\end{pmatrix}.
\end{equation*}
 We have $W_e(Q)\supset{\mathcal{D}}$  if and only if  $W_e(Q_{{\mathcal{H}}_{ns}})\supset{\mathcal{D}}$. By Lemma \ref{lemred} this is necessary. We may approximate $W_e(Q_{{\mathcal{H}}_{ns}})$  with sligthly larger essential numerical ranges, by using a positive  diagonalizable operator $R_{\varepsilon}$ such that $R_{\varepsilon} \ge R \ge R_{\varepsilon}-\varepsilon I$, for which
$$
W_e\left( \begin{pmatrix} I& 0 \\ R_{\varepsilon}&0\end{pmatrix} \right)=
W_e\left( \bigoplus_{n=1}^{\infty}\begin{pmatrix} 1& 0 \\ a_n&0\end{pmatrix}\right)
$$
where $\{a_n\}_{n=1}^{\infty}$ is a sequence of positive scalars, the eigenvalues of $R_{\varepsilon}$. By the previous step of the proof, this essential numerical range contains $\mathcal{D}$ if and only if $\overline{\lim}\, a_n \ge a_{\ast}$. If this holds for all $\varepsilon >0$, then $\| Q\|\ge \sqrt{5+2\sqrt{5}}$.  \end{proof}

\section{Unital, trace preserving positive linear maps}

\vskip 5pt\noindent
Unital positive linear maps $\Phi:\bM_n\to\bM_n$, the matrix algebra, which preserve the trace play an important role in matrix analysis and its applications. These maps are sometimes called doubly stochastic \cite{Ando}. 

We say that  $\Phi: {\mathrm{L}}({\mathcal{H}})\mapsto {\mathrm{L}}({\mathcal{H}})$ is trace preserving  if it preserves the trace ideal ${\mathcal{T}}$ and ${\mathrm{Tr}}\,\Phi(Z)= {\mathrm{Tr}}\, Z$ for all $Z\in{\mathcal{T}}$.

\vskip 5pt
\begin{cor}\label{cordoubly} Let $A\in {\mathrm{L}}({\mathcal{H}})$. The following two conditions are equivalent:
\begin{itemize}
\item[(i)] $W_e(A)\supset {\mathcal{D}}$.
\item[(ii)] For all  $X\in  {\mathrm{L}}({\mathcal{H}})$ with $\| X\|<1$, there exists a unital, trace preserving,  positive linear map $\Phi:  {\mathrm{L}}({\mathcal{H}})\to  {\mathrm{L}}({\mathcal{H}})$ such that $\Phi(A)=X$.
\end{itemize}
We may further require in (ii) that $\Phi$ is completely positive and {\rm{sot}}- and {\rm{wot}}-sequentially continuous.
\end{cor}

\vskip 5pt
\begin{proof} Assume (i). By Theorem \ref{pinching} we have a unitary $U: {\mathcal{H}}\to\oplus^{\infty}{\mathcal{H}}$ such that
$$
A\simeq UAU^* =
\begin{pmatrix} 
X &\ast & \cdots &\cdots  \\ 

\ast &X &\ast &\cdots  \\ 

 \vdots &\ast  & X &\ddots \\ 

\vdots  &\vdots  &\ddots  &\ddots \\ 
\end{pmatrix}.
$$
Now consider the map $\Psi: {\mathrm{L}}(\oplus^{\infty}{\mathcal{H}})\to {\mathrm{L}}({\mathcal{H}})$,
$$
\begin{pmatrix} 
Z_{1,1} & Z_{1,2} & \cdots  \\ 

Z_{2,1} &Z_{2,2} &\cdots   \\ 

 \vdots &\vdots &\ddots \\  
\end{pmatrix}
\mapsto
\sum_{i=1}^{\infty} 2^{-i} Z_{i,i}
$$
and define $\Phi: {\mathrm{L}}({\mathcal{H}})\to {\mathrm{L}}({\mathcal{H}})$ as $\Phi(T)=\Psi(UTU^*)$. Since both $\Psi$ and the unitary congruence with $U$ are sot- and wot-sequentially continuous, and trace preseverving, completely positive and unital, so is $\Phi$. Further $\Phi(A)=X$.

Assume (ii) and suppose that $z\notin W_e(A) $ and $|z|<1$ in order to reach a contradiction. If $z=|z| e^{i\theta}$, replacing $A$ by $ e^{-i\theta}A$, we may assume $1>z\ge 0$. Hence,
$$
W_e((A+A^*)/2)\subset (-\infty, z]
$$
and there exists a selfadjoint compact operator $L$ such that
$$
\frac{A+A^*}{2}  \le zI +L.
$$
This implies that $X:=\frac{1+z}{2}I$ cannot be in the range of $\Phi$ for any unital, trace preserving positive linear map. Indeed, we would have 
$$
\frac{1+z}{2}I= \frac{X+X^*}{2}=\Phi\left( \frac{A+A^*}{2}\right)\le zI +\Phi(L)
$$
which is not possible as $\Phi(L)$ is compact.
\end{proof}

\vskip 5pt
 In the finite dimensional setting, two Hermitian matrices $A$ and $X$ satisfy the relation $X=\Phi(A)$ for some positive, unital, trace preserving linear map if and only if $X$ is in the convex hull of the unitary orbit of $A$. In the infinite dimensional setting, if two Hermitian $A,X\in{\mathrm{L}}({\mathcal{H}})$ satisfy $W_e(A)\supset [-1,1]$ and $\| X\| \le 1$, then $X$ is in the norm closure of the unitary orbit of $A$. This is easily checked by approximating the operators with diagonal operators. Such an equivalence might not be brought out to the  setting of Corollary \ref{cordoubly}.

\vskip 5pt
\begin{question} Do there exist $A,X\in {\mathrm{L}}({\mathcal{H}})$ such that  $W_e(A)\supset {\mathcal{D}}$, $\|X\| <1$, and $X$ does not belong to the norm closure of the convex hull of the unitary orbit of $A$ ?
\end{question}

Here we mention a result of Wu \cite[Theorem 6.11]{Wu}: {\it If $A\in{\mathrm{L}}({\mathcal{H}})$ is not of the form scalar plus compact, 
then every $X\in{\mathrm{L}}({\mathcal{H}})$ is a linear combination of operators in the unitary orbit of $A$.}

 If one deletes the positivity assumption, the most regular class of linear maps on
${\mathrm{L}}({\mathcal{H}})$ might be given in the following definition.

\vskip 5pt
\begin{definition}\label{ultra} A linear map $\Psi:  {\mathrm{L}}({\mathcal{H}})\to  {\mathrm{L}}({\mathcal{H}})$ is said ultra-regular if it fulfills   two conditions:
\begin{itemize}
\item[(u1)] $\Psi(I)=I$ and $\Psi$ is trace preserving.
\item[(u2)] Whenever a sequence $A_n\to A$ for either the norm-, strong-,   or weak-topology, then we also have $\Psi(A_n)\to \Psi(A)$ for the same type of convergence.
\end{itemize}
\end{definition}

\vskip 5pt
Any ultra-regular linear map preserves the set of essentially scalar operators (of the form $\lambda I + K$ with $\lambda\in\bC$ and a compact operator $K$).  For its complement, we state our last corollary.

\vskip 5pt
\begin{cor}\label{corultra} Let $A\in {\mathrm{L}}({\mathcal{H}})$ be essentially nonscalar. Then, for all  $X\in  {\mathrm{L}}({\mathcal{H}})$ there exists a ultra-regular linear map $\Psi:  {\mathrm{L}}({\mathcal{H}})\to  {\mathrm{L}}({\mathcal{H}})$ such that $\Psi(A)=X$.
\end{cor}

\vskip 5pt
\begin{proof} An operator is essentially nonscalar precisely when its essential numerical range is not reduced to a single point. So, let $a,b\in W_e(A)$, $a\neq b$.
By a lemma of Anderson and Stampfli \cite{AS}, $A$ is unitarily equivalent to an operator on ${\mathcal{H}}\oplus{\mathcal{H}}$ of the form
$$
B=\begin{pmatrix} D & \ast \\
\ast & \ast
\end{pmatrix}
$$
where $D=\oplus_{n=1}^{\infty} D_n$, with two by two matrices $D_n$,
$$
D_n=
\begin{pmatrix} a_n & 0\\
0 & b_n\end{pmatrix}
$$
such that $a_n\to a$ and $b_n\to b$ as $n\to\infty$. We may assume that, for some $\alpha,\beta >0$, we have $\alpha >|a_n|+|b_n|$ and $|a_n-b_n|>\beta$. Hence there exist $\gamma >0$ and two by two intertible matrices $T_n$ such that, for all $n$,
$W(T_nD_nT_n^{-1})\supset {\mathcal{D}} $ and $\| T_n\| +\| T_n^{-1}\|\le \gamma$. So, letting $T=\left(\oplus_{n=1}^{\infty} T_n\right)\oplus I$, we obtain an invertible operator $T$ on ${\mathcal{H}}\oplus{\mathcal{H}}$ such that $W_e(TBT^{-1})\supset {\mathcal{D}}$.

Hence we have an invertible operator $S$ on ${\mathcal{H}}$ such that $W_e(SAS^{-1})\supset {\mathcal{D}}$. Therefore we may apply Corollary  \ref{cordoubly} and obtain a wot- and sot-sequentially continuous, unital, trace preserving map $\Phi$ such that $\Phi(SAS^{-1})=X$. Letting $\Psi(\cdot)=\Phi(S\cdot S^{-1})$ completes the proof.
\end{proof}

\vskip 5pt
We cannot find an alternative proof, not based on the pinching theorem, for Corollaries
\ref{cordoubly} and \ref{corultra}. 

 If we trust in Zorn, there exists a linear map $\Psi :{\mathrm{L}}({\mathcal{H}})\to  {\mathrm{L}}({\mathcal{H}})$ which satifies the condition (u1) but not the condition (u2). Indeed, let $\{a_p\}_{p\in\Omega}$ be a basis in the Calkin algebra ${\frak{C}}={\mathrm{L}}({\mathcal{H}})/{\mathrm{K}}({\mathcal{H}})$, indexed on an ordered set $\Omega$, whose first element $a_{p_0}$ is the image of $I$ by the canonical projection $\pi: {\mathrm{L}}({\mathcal{H}})\to {\frak{C}}$. Thus, for each operator $X$, we have a unique decomposition 
$
\pi(X)= \sum_{p\in\Omega} (\pi(X))_p a_p
$
with only finitely many nonzero terms. Further $(\pi(X))_{p_0} = 0$ if $X$ is compact, and $(\pi(I))_{p_0} = 1$.
We then define a map  $\psi:{\mathrm{L}}({\mathcal{H}})\to  {\mathrm{L}}({\mathcal{H}}\oplus {\mathcal{H}})$ by
$$
\psi(X)= \begin{pmatrix} X&0 \\ 0& (\pi(X))_{p_0}I\end{pmatrix}.
$$
Letting $\Psi(X)=V\psi(X)V^*$ where $V:{\mathcal{H}}\oplus {\mathcal{H}}\to {\mathcal{H}}$ is unitary, we obtain a linear map $\Psi :{\mathrm{L}}({\mathcal{H}})\to  {\mathrm{L}}({\mathcal{H}})$ which satifies  (u1) but not  (u2): it is not norm continuous.

 Let $\omega$ be a Banach limit on $l^{\infty}$ and define a map $\phi:l^{\infty}\to l^{\infty}$, $\{a_n\}\mapsto \{b_n\}$, where $b_1=\omega(\{a_n\})$ and $b_n=a_{n-1}$, $n\ge 2$. Letting $\Psi(X)=\phi(diag(X))$, where $diag(X)$ is the diagonal of $X\in {\mathcal{H}}$ in an orthonormal basis, we obtain a linear map $\Psi$ which is norm continuous, satisfies (u1) but not (u2): it is not strongly sequentially continuous.

However, it seems not possible to define explicitly a linear map $\Psi :{\mathrm{L}}({\mathcal{H}})\to  {\mathrm{L}}({\mathcal{H}})$ satisfying  (u1) but not  (u2).

\section{Pinchings in factors ?}

We discuss possible extensions to our results to a von Neumann algebra ${\mathfrak{R}}$ acting on a separable Hilbert space ${\mathcal{H}}$. First, we need to define an essential numerical range $W_e^{\mathfrak{R}}$ for  ${\mathfrak{R}}$.   Let $A\in{\mathfrak{R}}$. If ${\mathfrak{R}}$ is type-${\mathrm{III}}$, then $W_e^{\mathfrak{R}}(A):=W_e(A)$.  If ${\mathfrak{R}}$ is type-${\mathrm{II}}_{\infty}$, then
$$
W_e^{\mathfrak{R}}(A):=\bigcap_{K\in {\mathcal{T}}}  \overline{W}(A+K)
$$
where ${\mathcal{T}}$ is the trace ideal in ${\mathfrak{R}}$ (we may also use its norm closure ${\mathcal{K}}$, the "compact" operators in  ${\mathfrak{R}}$, or any dense sequence in ${\mathcal{K}}$)

\vskip 5pt
\begin{question} In Corollaries \ref{cor1seq}, \ref{cor2seq} and  \ref{idem2seq}, can we replace  ${\mathrm{L}}({\mathcal{H}})$  by a type-${\mathrm{II}}_{\infty}$ or  -${\mathrm{III}}$ factor ${\mathfrak{R}}$ with $W_e^{\mathfrak{R}}$ ? 
\end{question}

\vskip 5pt
\begin{question} In Corollaries \ref{corKS} and  \ref{corLW}, can we replace  ${\mathrm{L}}({\mathcal{H}})$  by a type-${\mathrm{II}}_{\infty}$ or  -${\mathrm{III}}$ factor ${\mathfrak{R}}$  with $W_e^{\mathfrak{R}}$ ? 
\end{question}

\vskip 5pt
\begin{question} In Corollaries  \ref{cordoubly} and  \ref{corultra},  can we replace  ${\mathrm{L}}({\mathcal{H}})$ by a type-${\mathrm{II}}_{\infty}$ factor ${\mathfrak{R}}$ with $W_e^{\mathfrak{R}}$ ?
\end{question}

\vskip 5pt
Recently, Dragan and Kaftal \cite{DK} obtained some decompositions for positive operators in von Neumann factors, which, in the case of ${\mathrm{L}}({\mathcal{H}})$ were first investigated in \cite{BL1}-\cite{BL2} by using Theorem \ref{pinching}.  This suggests that our questions dealing with a possible extension to  type-${\mathrm{II}}_{\infty}$ and  -${\mathrm{III}}$ factors also have an affirmative answer. In fact,  it seems pausible that Theorem \ref{pinching} and Theorem \ref{pinchingnormal} admit a version for such factors and this would  affirmatively answer  these questions.

Let ${\mathfrak{R}}$ be a type-${\mathrm{II}}_{\infty}$ or  -${\mathrm{III}}$ factor. 

\vskip 5pt
\begin{definition} A sequence $\{V_i\}_{i=1}^{\infty}$ of isometries in ${\mathfrak{R}}$ such that $\sum_{i=1}^{\infty} V_iV_i^*=I$ is called an isometric decomposition of ${\mathfrak{R}}$.
\end{definition}

\vskip 5pt
\begin{conjecture} Let $A\in{\mathfrak{R}} $ with $W_e^{\frak{R}}(A)\supset{\mathcal{D}}$ and  $\{X_i\}_{i=1}^{\infty}$  a sequence in ${\mathfrak{R}}$ such that
$\sup_i\|X_i\|<1$. Then, there exists an isometric decomposition  $\{V_i\}_{i=1}^{\infty}$ of ${\mathfrak{R}}$ such that
 $V_i^*AX_i=X_i$ for all $i$.
\end{conjecture}

\vskip 10pt

J.-C. Bourin,

Laboratoire de math\'ematiques,

Universit\'e de Franche-Comt\'e,

25 000 Besan\c{c}on, France.

jcbourin@univ-fcomte.fr

\vskip 10pt
Eun-Young Lee

 Department of mathematics,

Kyungpook National University,

 Daegu 702-701, Korea.

eylee89@knu.ac.kr

\end{document}